\documentclass{article}
\usepackage{graphicx}

\newtheorem{theorem}{Theorem}

\title{Deformations of crystal frameworks}
\author{Ciprian S. Borcea and Ileana Streinu}
\date{}

\begin{document}
\maketitle

\begin{abstract}
We apply our deformation theory of periodic bar-and-joint frameworks to tetrahedral crystal structures. The deformation space is investigated in detail for frameworks modelled on quartz, cristobalite and tridymite.
\end{abstract}


\section*{Introduction}

In this paper we present specific applications of our general deformation theory of periodic frameworks \cite{BS}.

\medskip
Considerations related to {\em framework flexibility} appear already in the early structural investigations based on $X$-ray crystallography \cite{G1, G2, P1, P2}. 
For framework materials, envisaged as corner sharing polyhedra, an intuitive notion of a `coordinated tilting' of the polyhedra is used 
in classifying similar structures \cite{Gla, M} or in studies of thermal and pressure effects. 
A most important area of theoretical and experimental studies where geometric models of deforming frameworks have been implicated is that concerned with
{\em displacive phase transitions} \cite{GD, Dol, D}.

\medskip
 Regarding the use of geometrical facts, it should be observed that, for most
framework structures, only a confined sample of geometrical possibilities has been explored
in the literature, typically one-parameter families which are intuitively `accessible'. The deformation theory developed in our paper \cite{BS} shows that one may expect, in general, a rich and diverse geometry. The present  undertaking describes the deformation spaces for tetrahedral periodic frameworks modeled on {\em quartz, cristobalite and tridymite}.

\section{The quartz framework}

\begin{figure}[h]
 \centering
 {\includegraphics[width=0.95\textwidth]{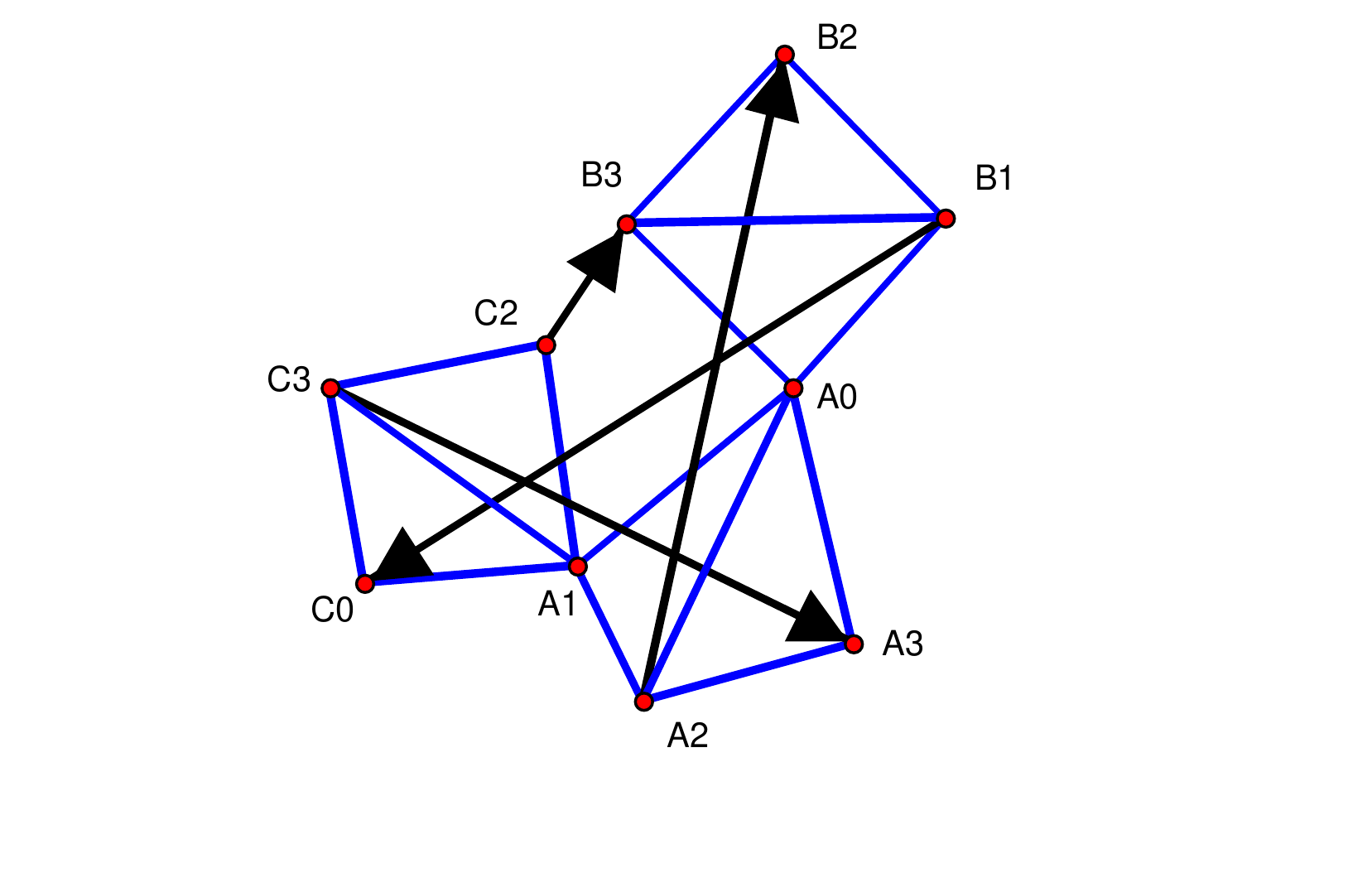}}
 \caption{A fragment of the tetrahedral framework of quartz. The periodicity lattice is generated by the four marked vectors, which must maintain a zero sum under deformation. The full framework is obtained by translating the depicted tetrahedra with all periods.}
 \label{FigQ}
\end{figure}

The ideal structure considered here is made of congruent regular tetrahedra. The oxygen atoms
would correspond with the vertices, each oxygen being shared by two tetrahedra. The silicon
atoms should be imagined at the centers of the tetrahedra. We shall examine all the geometric
deformations of the periodic framework described in Figure~\ref{FigQ}, without concern for
self-collision or any prohibition of a physical nature.

Equivalence under Euclidean motions is eliminated by assuming the tetrahedron marked $A_0A_1A_2A_3$ as fixed. Since all edges maintain their length, the positions of the two
tetrahedra which share the vertices $A_0$ and $A_1$ are completely described by two
orthogonal transformations $R_0$, respectively $R_1$ as follows: $R_0$ fixes $A_0$ and
takes $A_i$ to $B_i$, $i\neq 0$, while $R_1$ fixes $A_1$ and takes $A_j$ to $C_j$, $j\neq 1$.
The figure, by depicting only the `visible' edges, implies that both $R_0$ and $R_1$ are orientation reversing, that is, as orthogonal matrices $-R_0,-R_1\in SO(3)$.

\medskip \noindent
If we denote the edge vectors $A_i-A_0$ by $e_i, i=1,2,3$, we have:

$$ B_3-C_2=R_0e_3-(e_1+R_1(e_2-e_1)) $$

$$ A_3-C_3=e_3-(e_1+R_1(e_3-e_1)) $$

$$ B_2-A_2=R_0e_2-e_2 $$

$$ C_0-B_1=e_1-R_1e_1-R_0e_1 $$

\noindent
It follows that the dependency condition of a zero sum for these four generators of the periodicity lattice takes the form

\begin{equation}\label{0sum}
R_1(e_1-e_2-e_3)-R_0(e_1-e_2-e_3)=e_1+e_2-e_3
\end{equation}

\noindent
Under our regularity assumptions, the three vectors $R_1(e_1-e_2-e_3)$, $R_0(e_1-e_2-e_3)$
and $(e_1+e_2-e_3)$ have the same length and form an equilateral triangle. This restricts
$R_0(e_1-e_2-e_3)$ to the circle on the sphere of radius $||e_1-e_2-e_3||$ (which corresponde with an angle of $2\pi/3$ with $e_1+e_2-e_3$). Thus, $-R_0\in SO(3)$ is constrained to a
surface, which is differentiably a two-torus $(S^1)^2$. 

\medskip \noindent
For each choice of $-R_0$ on this torus, $R_1(e_1-e_2-e_3)$ is determined by (\ref{0sum}),
hence $-R_1$ is restricted to a circle $S^1$ in $SO(3)$. It follows that 

\begin{theorem}\label{quartz}
The deformation space of the ideal quartz framework is given by a three dimensional torus $(S^1)^3$ minus the degenerate cases when the span of the four vectors is less than three dimensional. 
\end{theorem}

\section{The cristobalite framework}

\begin{figure}[h]
 \centering
 {\includegraphics[width=0.70\textwidth]{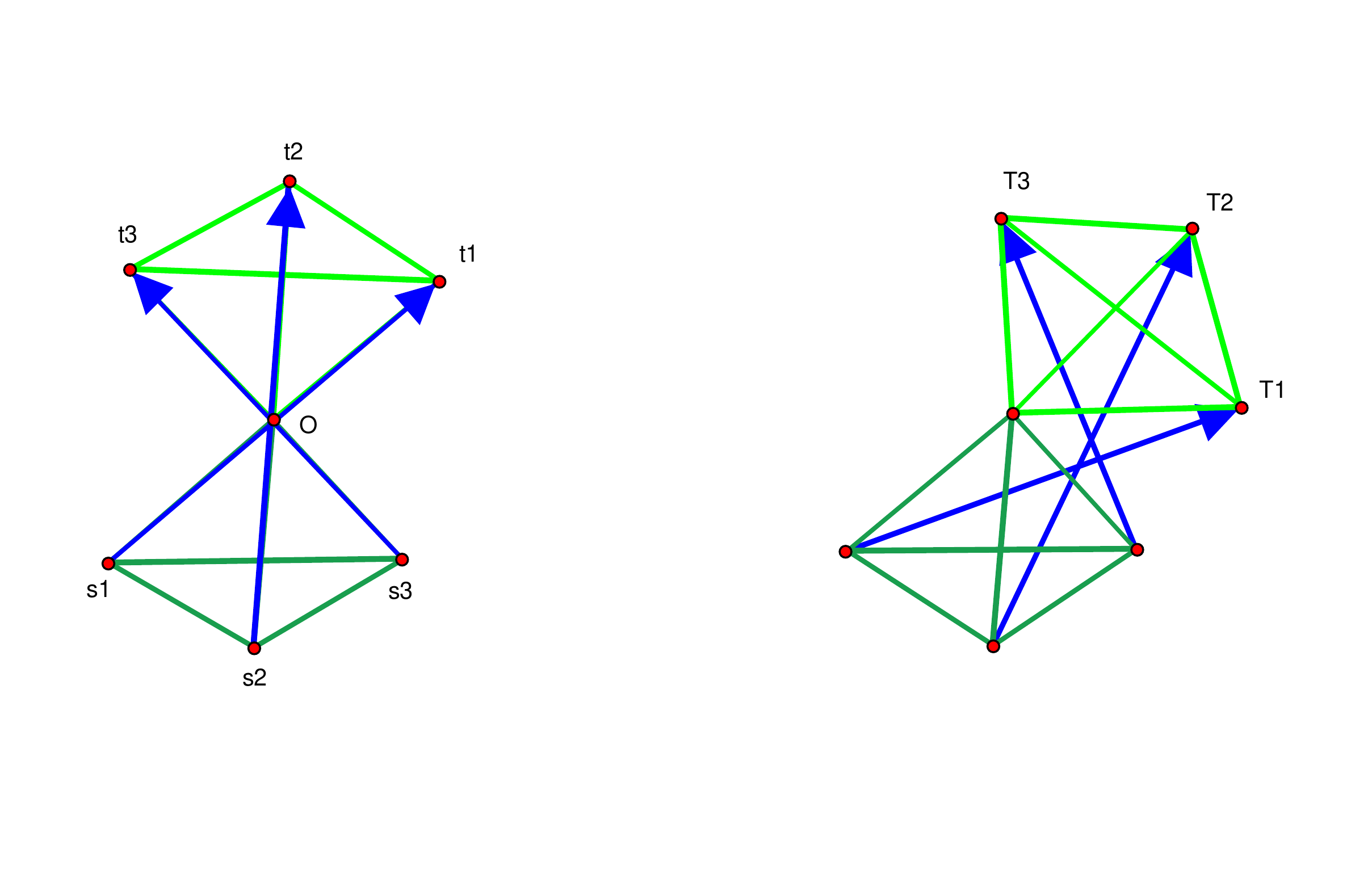}}
 \caption{Deforming the ideal cristobalite framework. The periodicity lattice is generated by the three vectors
 $\gamma_i=t_i-s_i$ which vary as the framework deforms.}
 \label{FigC}
\end{figure}

\medskip
The case of the `ideal $\beta$ cristobalite' structure illustrated in Figure~\ref{FigC} is already covered in \cite{BS}. The periodicity group of the framework is give by all the translational symmetries of the ideal crystal framework. As a result, there are $n=4$ orbits of vertices and
$m=12$ orbits of edges.

\begin{figure}[h]
 \centering
 {\includegraphics[width=0.60\textwidth]{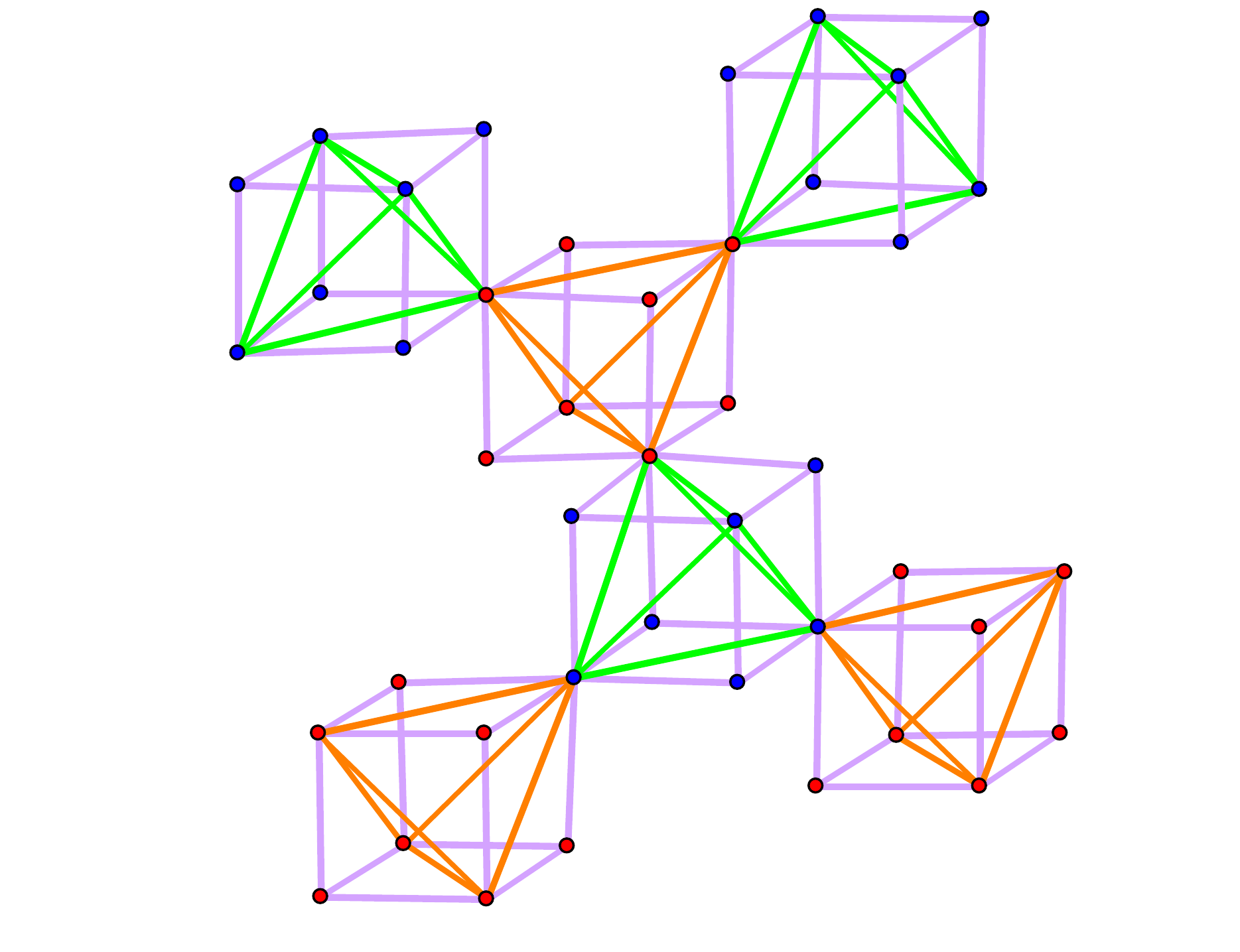}}
 \caption{ The ideal cristobalite framework (aristotype). The framework is made of vertex sharing regular tetrahedra.
 Cubes are traced only for suggestive purposes regarding symmetry and periodicity. See also Figure~\ref{FigC}.}
 \label{FigCC}
\end{figure}

Adopting the notations of Figure~\ref{FigCC}, we may assume the tetrahedron $Os_1s_2s_3$
as fixed and parametrize the possible positions of the other tetrahedon by a rotation around 
the origin $O$.  

\begin{theorem}\label{cristobalite}
The deformation space of the ideal high cristobalite framework is naturally parametrized by the open set in $SO(3)$ where the depicted generators remain linearly independent.
\end{theorem}

\section{The tridymite framework}

\begin{figure}[h]
 \centering
 {\includegraphics[width=0.80\textwidth]{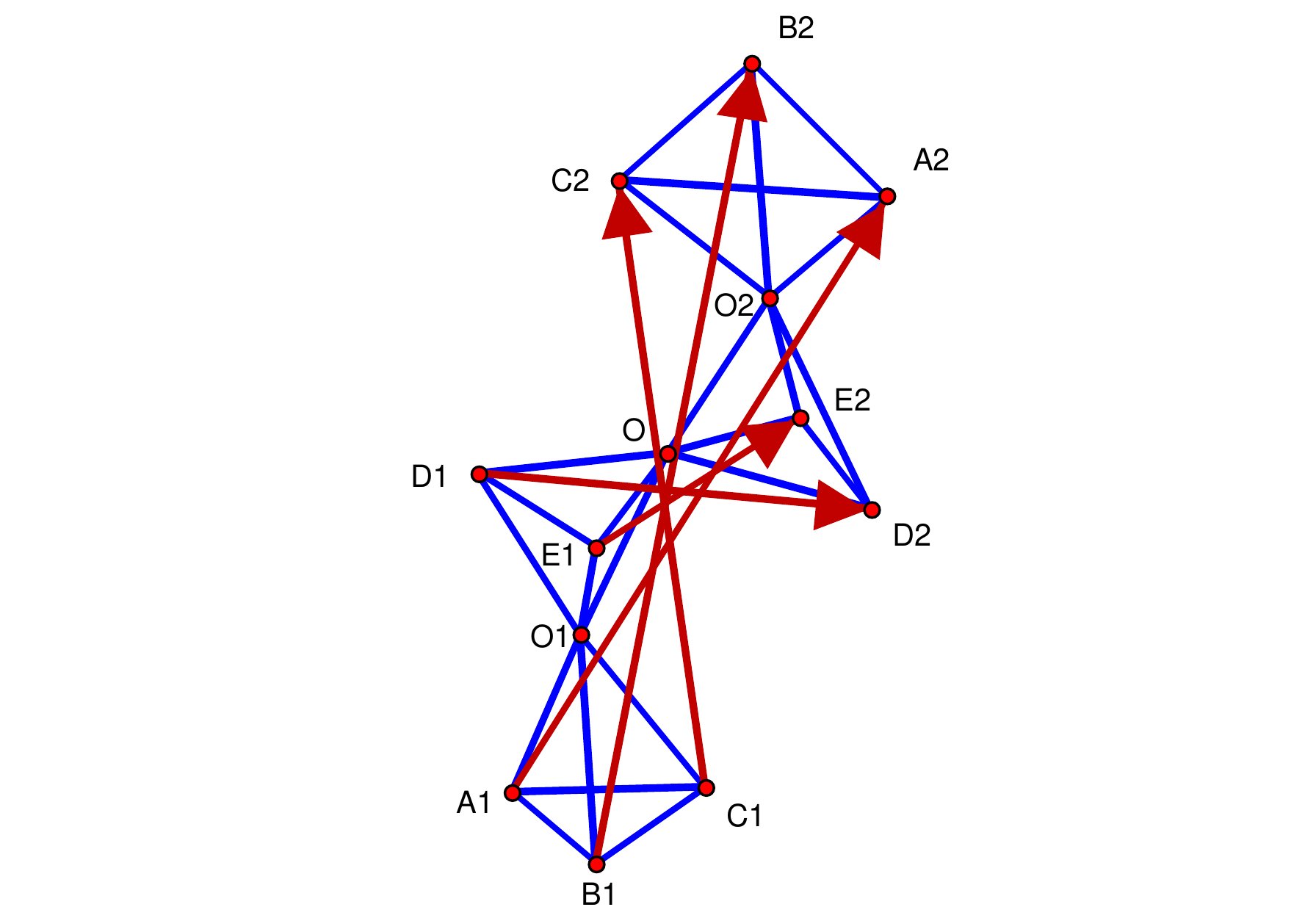}}
 \caption{The tetrahedral framework of tridymite. The periodicity lattice is generated by the marked
 vectors, subject to the relations
 $(C_2-C_1)+(D_2-D_1)=(A_2-A_1)$ and $(C_2-C_1)+(E_2-E_1)=(B_2-B_1)$. }
 \label{FigT}
\end{figure}

The tetrahedral framework $(G,\Gamma)$ of tridymite is depicted in Figure~\ref{FigT}. We consider the ideal case made of regular tetrahedra. The quotient graph has $|V/\Gamma|=8$
and $|E/\Gamma|=24$. All deformations can be described by three orthogonal transformations
(matrices) $R_0, R_1,R_2$ acting with centers at $O, O1$ and respectively $O2$.  With $O$ as
the origin and the tetrahedron $OD_1E_1O_1$ assumed fixed, we put:

$$ O_1=f_0, \ \ D_1=f_1 \ \ \mbox{and} \ \ f_1=f_2 $$

\noindent
Then, our orthogonal transformations are determined by the following relations:

$$ O_2=R_0f_0,\ \ D_2=R_0f_1 \ \ \mbox{and} \ \ f_2=R_0f_2 $$
$$ A_1=f_0+R_1(f_1-f_0),\ \ B_1=f_0+R_1(f_2-f_0), \ \ C_1=f_0-R_1f_0 $$
$$ A_2=R_0f_0+R_2R_0(f_1-f_0),\ B_2=R_0f_0+R_2R_0(f_2-f_0),\ C_2=R_0f_0-R_2R_0f_0 $$

\noindent
As a result, the two linear dependence relations between the six depicted periods take the form:

\begin{equation}\label{Tdepend}
(I-R_0-R_1+R_2R_0)f_i=0, \ \ i=1,2
\end{equation}

\noindent
where $I$ denotes the identity. We note that the ideal high tridymite structure
(the {\em aristotype}) corresponds to $R_0=-I$ and $R_1=R_2$ the reflection in the
plane $span(f_1,f_2))$.
 
\medskip
We shall describe the deformation space in a neighbourhood of this high tridymite structure.
We put $-R_0=Q$, $R_1=Q_1$ and $-R_2R_0=Q_2$, so that (\ref{Tdepend}) becomes

\begin{equation}\label{TdependBis}
 I+Q=Q_1+Q_2 \ \ \mbox{on}\ \ span(f_1,f_2)
\end{equation}

\noindent
with $Q,-Q_1,-Q_2\in SO(3)$. Since  the orthogonal transformations $Q,Q_1,Q_2$ are
completely determined by their values on two vectors $e_1,e_2$ of a Cartesian frame
with $span(e_1,e_2)=span(f_1,f_2)$, we have to solve the system

\begin{equation}\label{Tsystem}
 e_i+Qe_i=Q_1e_i+Q_2e_i \ \  i=1,2
\end{equation}

\noindent
where we assume $Q\in SO(3)$ given in a neighbourhood of the identity, and look for
solutions $Q_1,Q_2$.

\begin{figure}[h]
 \centering
 {\includegraphics[width=0.80\textwidth]{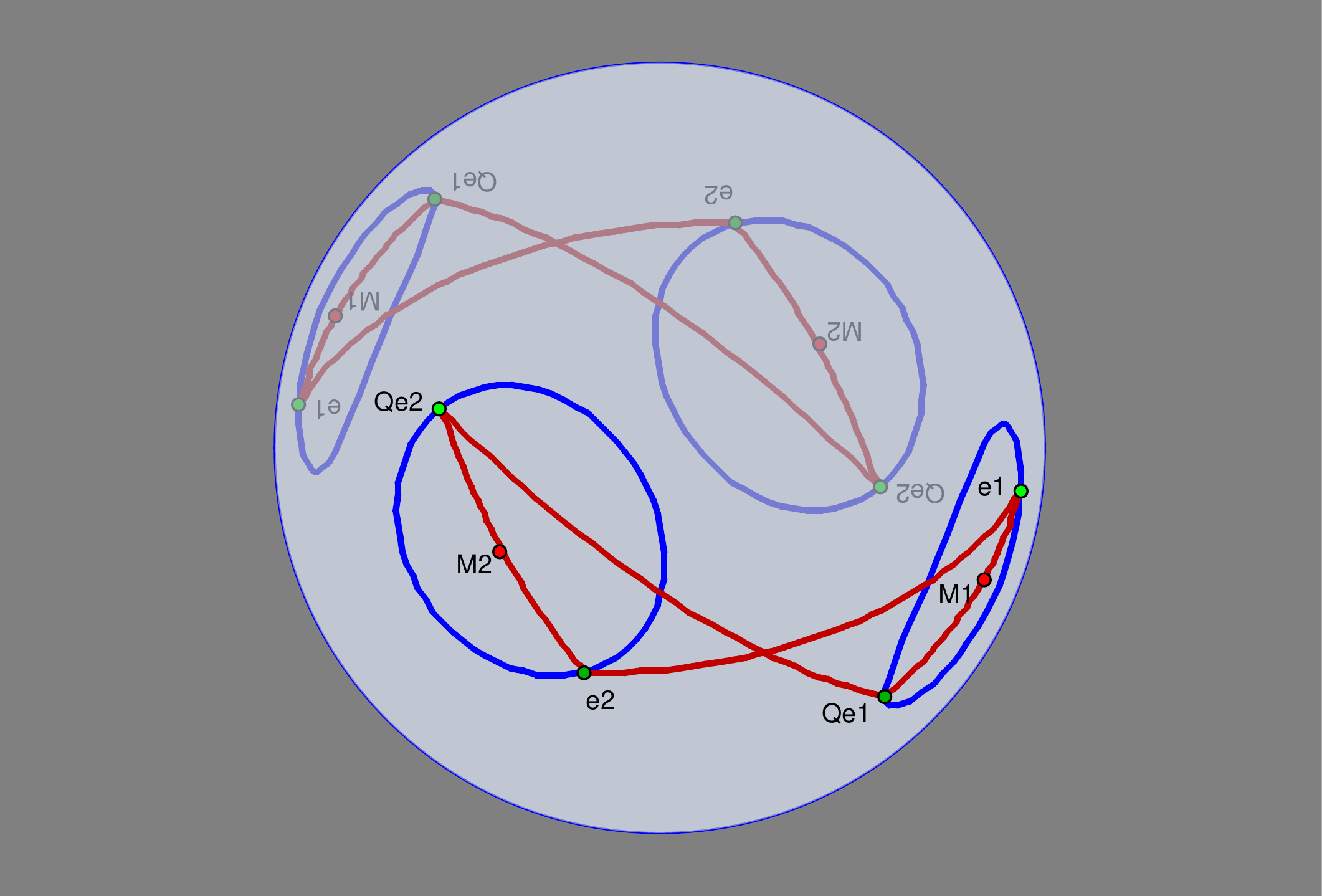}}
 \caption{The  spherical four-bar mechanism associated to the system (\ref{Tsystem}). }
 \label{FigS}
\end{figure}

\medskip
 We may interpret this system as a problem about a {\em spherical four-bar mechanism} in the
following way. All the vectors implicated in (\ref{Tsystem}) are unit vectors and can be depicted as points on the unit sphere $S^2$. For a given $Q$, we mark by $M_i$ the midpoint
of the spherical geodesic segment $[e_i,Qe_i]$ and trace the circle with center $M_i$ and diameter $[e_i,Qe_i]$. This is illustrated in Figure~\ref{FigS}.

\medskip \noindent
It is an elementary observation that any solution $Q_1e_i$ and $Q_2e_i$
determines diameters of the corresponding  circles for $i=1,2$, with the two geodesic arcs
$[Q_ke_1,Q_ke_2]$, like $[e_1,e_2]$ and $[Qe_1,Qe_2]$, of length $\pi/2$. Thus, the two spherical quadrilaterals with vertices at $e_1, Qe_1, Qe_2, e_2$ and respectively $Q_1e_1$, $Q_2e_1$, $Q_2e_2$, $Q_1e_2$  are two configurations of the same four-bar mechanism and
moreover, the distance between the midpoints of the opposite edges represented by
diameters is the same.

\medskip \noindent
It follows from the theory of the spherical four-bar mechanism that, for a generic $Q$ near the identity of $SO(3)$, the abstract configuration space is made of two loops which correspond by reflecting the corresponding realizations. Each loop component has two configurations with
the prescibed distance $[M_1M_2]$. Thus, there are four configurations with the prescribed
distance. 

\medskip \noindent
We observe that if we replace $Q_1$ by $Q_2$ and $Q_2$ by $Q_1$ in the labeling of the
vertices  of a realization, the orientation is reversed, hence the configuration belongs to the 
other component. Thus, the two obvious solutions of (\ref{Tsystem}), namely

$$ Q_1e_i=e_i,\  Q_2e_i=Qe_i \ \ \ \mbox{and} \ \ \  Q_1e_i=Qe_i, \ Q_2e_i=e_i , \ \ \ i=1,2 $$

\noindent
correspond with configurations on the two different loop components, as do the 
remaining two, which are also paired by relabeling. This discussion shows that
all four solutions are obtained from the quadrilateral $e_1, Qe_1, Qe_2, e_2$ and
its reflection in  the geodesic $[M_1,M_2]$, by the two relabelings with  $Q_1$ and $Q_2$
possible in each case.

\begin{figure}[h]
 \centering
 {\includegraphics[width=0.80\textwidth]{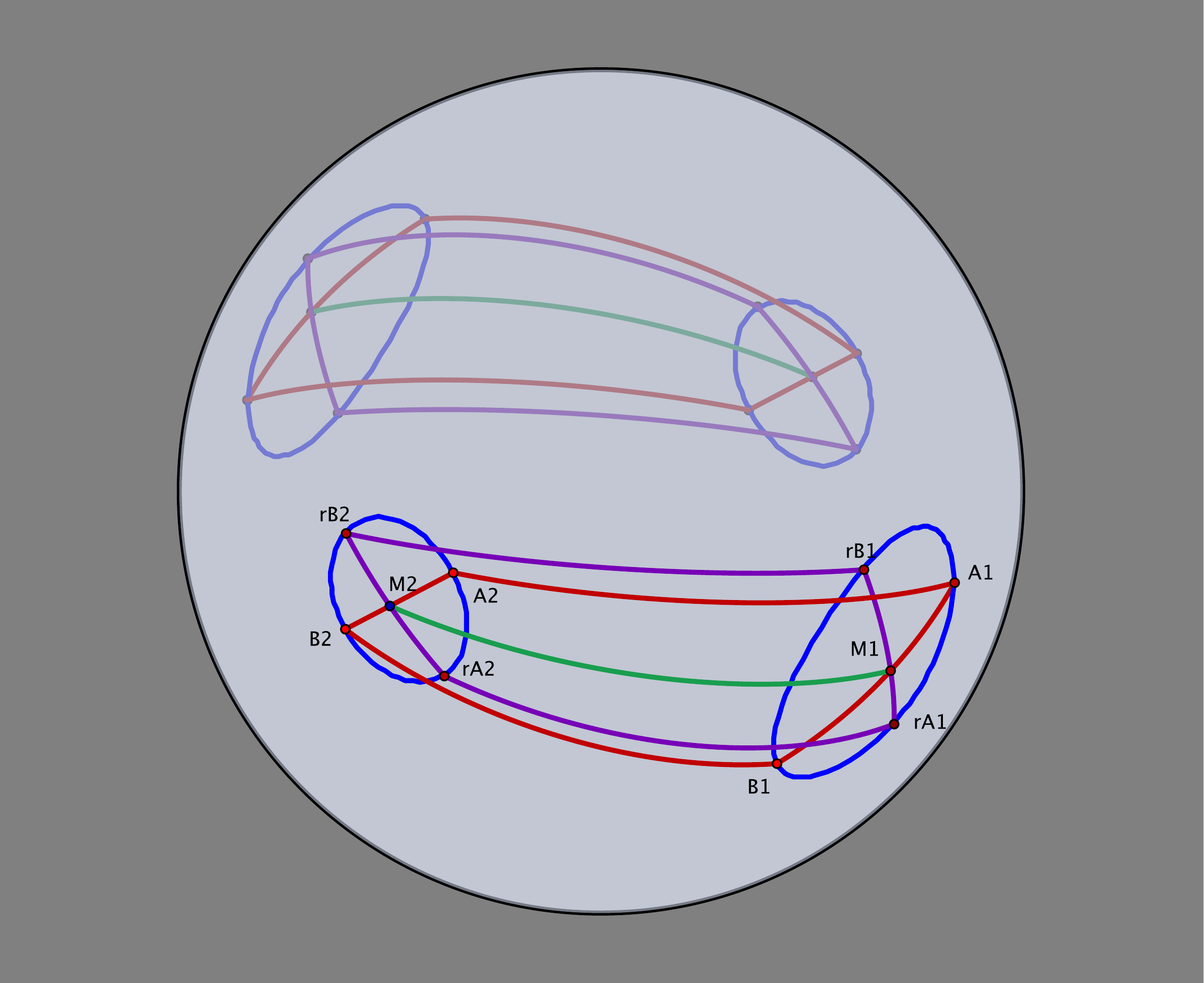}}
 \caption{Spherical  four-bar mechanism and reflection in $[M_1,M_2]$. }
 \label{FigSR}
\end{figure}

\medskip \noindent
In Figure~\ref{FigSR} we have depicted the quadrilateral $e_1, Qe_1, Qe_2, e_2$ as 
$A_1B_1B_2A_2$, with reflection in $[M_1M_2]$ marked as $rA_1,rB_1,rB_2,rA_2$.
Then,, the solutions $(Q_1e_1, Q_1e_2,Q_2e_1,Q_2e_2)$ of the system (\ref{Tsystem})
are the following four solutions: 
$(A_1,A_2,B_1,B_2)$, $(B_1,B_2,A_1,A_2)$,  $(rA_1,rA_2,rB_1,rB_2)$ and 
$(rB_1,rB_2,rA_1,rA_2)$.

\medskip \noindent
We may summarize our result as follows.

\begin{theorem}
The deformation space of the tridymite framework is singular in a neighbourhood of the
aristotype and can be represented as a ramified covering with four sheets of a three-dimensional
domain. There is a natural $Z_2\times Z_2$ action on this covering which fixes the aristotype framework.
\end{theorem}

\medskip \noindent
Indeed,  the two involutions, inverting the labeling and reflecting in $[M_1,M_2]$, commute
ang give a $Z_2\times Z_2.$ action on the covering. The dimension of the tangent space
at the aristotype framework is computed from the linear version of (\ref{Tsystem}) and is
six.

\vspace{0.2in}
 
   Ciprian S. Borcea\\
   Department of Mathematics\\
   Rider University\\
   Lawrenceville, NJ 08648, USA

\vspace{0.2in} 

   Ileana Streinu\\
   Department of Computer Science\\
   Smith College\\ 
   Northampton, MA 01063, USA

\end{document}